\documentclass[letterpaper, 10 pt, conference]{ieeeconf}  %

\IEEEoverridecommandlockouts                              %
\usepackage{amsmath}
\usepackage{amssymb}
\usepackage{amsfonts}
\usepackage{graphicx}
\usepackage[hidelinks]{hyperref}
\usepackage{mathtools}
\usepackage{booktabs}
\usepackage{algorithm}
\usepackage{algpseudocode}
\usepackage{float}
\usepackage{xcolor}
\usepackage{comment}

\DeclareMathOperator{\chol}{chol}

\newcommand{\rev}[1]{#1}
\newenvironment{revenv}{\bgroup}{\egroup}

\title{\LARGE \bf
Exploiting Multistage Optimization Structure in Proximal Solvers
}

\author{Roland Schwan, Daniel Kuhn, and Colin N. Jones
\thanks{Open source: \url{https://github.com/PREDICT-EPFL/piqp}}%
\thanks{Claude.ai has been used to improve the syntax and grammar of several paragraphs in the manuscript. This work was supported as a part of NCCR Automation, a National Centre of Competence in Research, funded by the Swiss National Science Foundation (grant number 51NF40\_225155). 
Roland Schwan, and Colin N. Jones are with the Automatic Control Lab, EPFL, Switzerland. Roland Schwan and Daniel Kuhn are with the Risk Analytics and Optimization Chair, EPFL, Switzerland. \tt{\{roland.schwan, daniel.kuhn, colin.jones\}@epfl.ch}}
}%

\begin{document}

\maketitle
\thispagestyle{empty}
\pagestyle{empty}

\begin{abstract}
This paper presents an efficient structure-exploiting algorithm for multistage optimization problems. The proposed method extends existing approaches by supporting full coupling between stages and global decision variables in the cost, as well as equality and inequality constraints. \rev{The algorithm is implemented as a new backend in the PIQP solver and leverages} a specialized block-tri-diagonal-arrow Cholesky factorization within a proximal interior-point framework to handle the underlying problem structure efficiently. The implementation features automatic structure detection and seamless integration with existing interfaces. Numerical experiments demonstrate significant performance improvements, achieving up to 13x speed-up compared to a generic sparse backend and matching/exceeding the performance of the state-of-the-art specialized solver HPIPM. The solver is particularly effective for applications such as model predictive control, robust scenario optimization, and periodic optimization problems.
\end{abstract}

\section{Introduction}

The efficient solution of convex optimization problems is crucial across engineering disciplines, particularly in real-time applications. An important category is the class of multistage optimization problems, which arise naturally in various domains, including model predictive control (MPC) and multistage portfolio optimization. These problems are characterized by a sequence of interconnected stages, often coupled through both the constraints and the cost function.

\rev{This paper presents an approach for solving general multistage optimization problems, supporting full coupling between stages and global decision variables in the cost and all constraint types.}
While numerous solvers have been developed for multistage optimization problems, most impose restrictive assumptions on the underlying problem structure. A comprehensive list of multistage quadratic programming (QP) solvers can be found in Table~\ref{tbl:qp_comparison}. These solvers typically support coupling between consecutive stages only through equality constraints or forward system dynamics, limiting their applicability to more general problems. This limitation extends to non-linear optimization solvers like FATROP \cite{vanroye2023} or FORCES NLP \cite{zanelli2020}, \rev{and to tree-sparse optimization problems \cite{steinbach2003} and their specialized solvers \cite{huebner2017,frison2020}.} Though \cite{pacaud2024} and \cite{shin2024} consider global decision variables, they assume independent stages without coupling. Many problems, including periodic optimization \cite{limon2016}, robust scenario MPC \cite{lucia2013}, or tracking MPC \cite{limon2008}, can be naturally formulated as multistage optimization problems with global decision variables.

The main contributions \rev{of this paper} are: (1) the exploitation of the block-tri-diagonal-arrow structure of the KKT matrix that naturally arises from this coupling and (2) the automatic detection of this structure from standard sparse problem formulations. Although we present our method within a proximal interior-point framework, the structure-exploiting factorization can be easily integrated into various proximal-based solvers, \rev{such as} OSQP \cite{stellato2020}. Additionally, we provide seamless integration with existing interfaces through automatic structure detection that is invariant to constraint ordering, significantly reducing implementation burden.

The proposed method is implemented as a new backend in the PIQP solver \cite{schwan2023}, combining the robustness of proximal interior-point methods with the efficiency of structure-exploiting linear algebra. Our numerical experiments demonstrate significant performance improvements, achieving up to 13x speed-up compared to the generic sparse backend. Notably, the solver maintains or exceeds the performance of state-of-the-art solvers like HPIPM \cite{frison2020}.

Section~\ref{sec:problem_formulation} presents the problem formulation. Section~\ref{sec:prximal_ip} reviews the proximal interior-point method underlying our approach. Section~\ref{sec:solving_kkt} introduces our structure-exploiting factorization routines, followed by implementation details in Section~\ref{sec:numerical_implementation}. Finally, numerical results are presented in Section~\ref{sec:numerical_examples}, showcasing the solver's performance.

\subsection*{Notation}

We denote the set of $n$-dimensional real vectors by $\mathbb{R}^n$, and $n \times m$-dimensional real matrices by $\mathbb{R}^{n \times m}$. Moreover, the subspace of symmetric matrices in $\mathbb{R}^{n \times n}$ is denoted by $\mathbb{S}^n$, and the cone of positive semi-definite matrices in $\mathbb{S}^n$ by $\mathbb{S}_+^n$. The identity matrix is denoted by $I_n \in \mathbb{R}^{n \times n}$ and a vector filled with ones by $\mathbf{1}_n \in \mathbb{R}^n$.

\section{Problem Formulation} \label{sec:problem_formulation}

We consider the multistage optimization problem
\begin{equation} \label{eq:stage_problem}
\begin{aligned}
    \min_{x, g} \quad & \sum_{i=0}^{N-1} \ell_i(x_i,x_{i+1},g) + \ell_N(x_N,g) \\
    \text {s.t.}\quad & A_ix_i+B_ix_{i+1}+E_ig = b_i,\hspace{0.5em} i=0,\dots,N-1, \\
    & C_ix_i+D_ix_{i+1}+F_ig \leq h_i,\hspace{0.40em} i=0,\dots,N-1, \\
    & A_Nx_N + E_Ng = b_N, \\
    & D_Nx_N + F_Ng \leq h_N,
\end{aligned}
\end{equation}
with $N$ stages and coupled stage cost
\begin{table*}[htbp]
    \centering
    \caption{Comparison of Structure-Exploiting Multistage QP Solvers} \label{tbl:qp_comparison}
    \vspace{-1em}
    \begin{tabular}{lccccc}
        \toprule
        Characteristics & PIQP & Forces \cite{domahidi2012} & HPIPM \cite{frison2020} & QPALM-OCP \cite{lowenstein2024} & ASIPM \cite{frey2020} \\
        \midrule
        Algorithm & PALM-IP & IP & IP & PALM & Inexact IP \\
        Stage-wise Structure & Yes & Yes & Yes & Yes & Yes \\
        Global Variables & Yes & No & No & No & No \\
        Inter-stage Coupling & Full & Equality only & Forward dynamics & Forward dynamics & Forward dynamics \\
        Implementation & C++ & C & C & Matlab/Codegen & C \\
        License Type & BSD-2 & Commercial & BSD-2 & Closed & Closed \\
        \bottomrule
    \end{tabular}
    \vspace{-1em}
    \label{tab:solver_comparison}
\end{table*}
\begin{equation*}
    \ell_i(x_i,x_{i+1},g)\! \coloneqq\! \frac{1}{2}\!\begin{bmatrix}
x_i \\
x_{i+1} \\
g
\end{bmatrix}^\top\!\!\begin{bmatrix}
Q_i & S_i^\top & T_i^\top \\
S_i & 0 & 0 \\
T_i & 0 & 0
\end{bmatrix}\!\!\!\begin{bmatrix}
x_i \\
x_{i+1} \\
g
\end{bmatrix} \!+ c_i^\top x_i,
\end{equation*}
and terminal cost
\begin{equation*}
    \ell_N(x_N,g) \coloneqq \frac{1}{2}\begin{bmatrix}
x_N \\
g
\end{bmatrix}^\top\begin{bmatrix}
Q_N & T_N^\top \\
T_N & Q_g
\end{bmatrix}\begin{bmatrix}
x_N \\
g
\end{bmatrix}+ c_N^\top x_N + c_g^\top g,
\end{equation*}
where $x_i\in\mathbb{R}^{n_i}$ are the stage-wise decision variables, $g\in\mathbb{R}^{n_g}$ is a global decision variable, and $N \in \mathbb{N}$ is the horizon. The matrices $Q_i\in\mathbb{S}_+^{n_i}$, $S_i\in\mathbb{R}_+^{n_{i+1}\times n_i}$, and $T_i\in\mathbb{R}_+^{n_g\times n_i}$ together with $c_i\in\mathbb{R}^{n_i}$ and $c_g\in\mathbb{R}^{n_g}$ encode the coupled cost. The stage-wise variables are also coupled through equality and inequality constraints encoded by $A_i \in \mathbb{R}^{p_i\times n_i}$, $B_i \in \mathbb{R}^{p_i\times n_{i+1}}$, $E_i \in \mathbb{R}^{p_i\times n_g}$, $b_i \in \mathbb{R}^{p}$, $C_i \in \mathbb{R}^{m_i\times n_i}$, $D_i \in \mathbb{R}^{m_i\times n_{i+1}}$, and $h_i \in \mathbb{R}^{m_i\times n_g}$, respectively. Note that the formulation \eqref{eq:stage_problem} is more general than the typical optimal control problem (OCP) formulations in the literature, as the stages are not only coupled through the system dynamics. \rev{Note that $x_i$ in \eqref{eq:stage_problem} differs from the usual OCP formulations as it would contain both the state and input in this case.}

We can reformulate problem \eqref{eq:stage_problem} into a general QP form
\begin{equation} \label{eq:general_qp}
\begin{aligned}
\min_{x,s} \quad & \frac{1}{2} x^\top P x + c^\top x \\
\text{s.t.}\quad & Ax=b, \\
& Gx - h + s = 0, \\
& s \geq 0,
\end{aligned}
\end{equation}
where the stages $x_i$ and the global variable $g$ are collected into the decision variable $x=(x_0,\dots,x_N,g)\in\mathbb{R}^n$ with $n=\sum_{i=0}^{N}n_{x_i}+n_g$. The cost is now defined by a block-tri-diagonal-arrow matrix $P\in\mathbb{R}^{n \times n}$ and $c\in\mathbb{R}^n$. By introducing the slack variable $s\in\mathbb{R}^m$ with $m=\sum_{i=0}^{N}m_{i}$, we lift the affine inequality constraints to equality constraints, leaving us with a simple non-negativity constraint and the equality constraints defined by $A\in\mathbb{R}^{p \times n}$, $b\in\mathbb{R}^{p}$, $G\in\mathbb{R}^{m \times n}$, and $h\in\mathbb{R}^{m}$, with $p=\sum_{i=0}^{N}p_{i}$.

\section{Proximal Interior-Point} \label{sec:prximal_ip}

This section introduces the interior-point method that was developed in \cite{schwan2023}. While the factorization algorithms presented in this paper should generalize to most proximal-based solvers, we demonstrate the effectiveness of the proposed method using the algorithm in PIQP.

The proximal augmented Lagrangian of problem \eqref{eq:general_qp} is given as
\begin{equation}
\begin{aligned}
\mathcal{L}_{\rho, \delta}(x,s&; \xi,\lambda,\nu) \coloneqq \frac{1}{2} x^\top P x + c^\top x + \frac{\rho}{2}\left\|x-\xi\right\|_2^2\\
&+\lambda^\top (A x-b)+\frac{1}{2 \delta}\|A x-b\|_2^2 \\
&+\nu^\top (Gx-h+s)+\frac{1}{2 \delta}\|Gx-h+s\|_2^2,
\end{aligned}
\end{equation}
where the variables $\lambda \in \mathbb{R}^p$ and $\nu \in \mathbb{R}^m$ are the Lagrange multipliers of the equality constraints in \eqref{eq:general_qp}, $\xi \in \mathbb{R}^n$ acts as a regularizer, and $\rho, \delta > 0$ are penalty parameters. The proximal method of multipliers, proposed in \cite{rockafellar1976}, is then given by
\begin{subequations}
\begin{align}
(\xi^{k+1},s^{k+1}) & \in \underset{x,s\geq 0}{\text{argmin}}\;\mathcal{L}_{\rho_k, \delta_k}(x,s; \xi^k,\lambda^k,\nu^k),  \label{eq:argmin_pmm} \\
\lambda^{k+1} & = \lambda^k+\frac{1}{\delta_k}\left(A \xi^{k+1}-b\right), \label{eq:lam_update} \\
\nu^{k+1} & = \nu^k+\frac{1}{\delta_k}\left(G \xi^{k+1}-h+s^{k+1}\right), \label{eq:nu_update}
\end{align}
\end{subequations}
where $\xi$ is used to represent the $x$ variable in the proximal method of multipliers.

PIQP follows the method \rev{outlined} in \cite{pougkakiotis2021}, where one iteration of the interior-point method is applied to \eqref{eq:argmin_pmm}. More specifically, we replace the inequality constraint in \eqref{eq:argmin_pmm} with the log-barrier $\Phi_\mu(s) = -\mu \sum_{i=1}^m \ln [s]_i$ resulting in the minimization problem
\begin{equation} \label{eq:min_pmm_unconst}
\underset{x,s}{\text{min}}\;\mathcal{L}_{\rho_k, \delta_k}(x,s; \xi^k,\lambda^k,\nu^k)+\sigma_k\cdot \Phi_{\mu_k}(s),
\end{equation}
where $[s]_i$ denotes the $i$-th element of $s$, $\mu>0$ is the barrier parameter, and $\sigma_k \in (0,1]$ is a scaling parameter. The resulting first-order optimality conditions of \eqref{eq:min_pmm_unconst} are then given by
\begin{subequations} \label{eq:kkt_conditions}
\begin{align}
Px + c + \rho_k(x - \xi^k) + A^\top y + G^\top z & =0, \\
Ax - \delta_k (y-\lambda^k) - b & =0, \\
Gx - \delta_k (z-\nu^k) - h + s & =0, \\
s \circ z - \sigma_k\mu_k \mathbf{1}_m & =0, 
\end{align}
\end{subequations}
with auxiliary variables $y \coloneqq \lambda^k+\frac{1}{\delta_k}(A x-b)$ and $z \coloneqq \nu^k+\frac{1}{\delta_k}(Gx - h + s)$, simplifying the equations and the updates \eqref{eq:lam_update} and \eqref{eq:nu_update}. Here, the $\circ$ operator indicates elementwise multiplication. Solving \eqref{eq:kkt_conditions} using Newton's method \rev{with exact line search} amounts to solving the following linear equation at iteration $k$
\begin{equation} \label{eq:kkt_full}
\begin{bmatrix}    
P + \rho_k I_n & A^\top & G^\top & 0 \\
A & -\delta_k I_p & 0 & 0 \\
G & 0 & -\delta_k I_m & I_m \\
0 & 0 & S^k & Z^k
\end{bmatrix}
\begin{bmatrix}
\Delta x^k \\
\Delta y^k \\
\Delta z^k \\
\Delta s^k
\end{bmatrix}
=
\begin{bmatrix}
r^k_x \\
r^k_y \\
r^k_z \\
r^k_s
\end{bmatrix}
\end{equation}
with initialization $(x^k,y^k,z^k,s^k)$, $S^k, Z^k \in \mathbb{R}^{m \times m}$ are diagonal matrices with $s^k$, $z^k$ on its diagonal, and residuals
\begin{equation*}
\begin{aligned}
r^k_x &= - (Px^k + c + \rho_k(x^k-\xi^k) + A^\top y^k + G^\top z^k), \\
r^k_y &= - (Ax^k + \delta_k(\lambda^k-y^k) - b), \\
r^k_z &= - (Gx^k + \delta_k(\nu^k - z^k) - h + s^k), \\
r^k_s &= -s^k \circ z^k + \sigma_k \mu_k\mathbf{1}_m.
\end{aligned}
\end{equation*}

Eliminating $\Delta s^k$ in \eqref{eq:kkt_full} results in a reduced system of equations
\begin{equation*}
\underbrace{\begin{bmatrix}
P + \rho^k I_n &\!\!\! A^\top &\!\!\! G^\top \\
A &\!\!\! -\delta^k I_p &\!\!\! 0 \\
G &\!\!\! 0 &\!\!\! -(W^k + \delta^k I_m)
\end{bmatrix}}_{J(s^k,z^k)}
\begin{bmatrix}
\Delta x^k \\
\Delta y^k \\
\Delta z^k
\end{bmatrix}
=\begin{bmatrix}
r^k_x \\
r^k_y \\
\bar{r}^k_z
\end{bmatrix}
\end{equation*}
with the Nesterov-Todd scaling $W^k = (Z^k)^{-1} S^k$ and $\bar{r}^k_z = r^k_z - (Z^k)^{-1} r^s_k$ \cite{nocedal2006}. The slack direction $\Delta s^k$ can then be reconstructed with $\Delta s^k = (Z^k)^{-1}(r^k_s - S^k \Delta z^k)$.

In the remainder of the paper, we will focus on how to solve $J(s^k,z^k)$ as efficiently as possible, exploiting the underlying structure. Thus, we will not go into further details concerning the interior-point method itself. The interested reader is referred to \cite{schwan2023} for more details.

Note that the structure of $J(s^k,z^k)$ is ubiquitous and is present in almost all proximal-based solvers. For example, OSQP \cite{stellato2020} has the exact same KKT structure. Hence, the factorization scheme presented in the next section is directly transferable.

\section{Solution of the Linearized KKT System} \label{sec:solving_kkt}

The symmetry of $J(s^k,z^k)$ enables the use of a sparse LDL factorization, as the regularization terms ensure the matrix is quasi-definite \cite{vanderbei1995}. While this approach was adopted in the original PIQP solver due to its general treatment of $P$, $A$, and $G$ without exploiting their structure, the multistage nature of our problem presents an opportunity for more efficient computation.

Eliminating $\Delta z^k$ gives the so-called \textit{augmented system}
\begin{equation*}
\begin{bmatrix}
\Phi &\!\!\! A^\top \\
A &\!\!\! -\delta^k I_p
\end{bmatrix}
\begin{bmatrix}
\Delta x^k \\
\Delta y^k
\end{bmatrix}
=\begin{bmatrix}
\bar{r}^k_x \\
r^k_y
\end{bmatrix}
\end{equation*}
where
\begin{equation}
\begin{aligned}
    \Phi &= P + \rho^k I_n + G^\top (W^k+\delta_kI_m)^{-1} G, \\
    \bar{r}^k_x &= r^k_x + G^\top (W^k+\delta_kI)^{-1} \bar{r}^k_z, \\
    \Delta z^k &= (W^k+\delta_kI)^{-1}(G \Delta x^k - \bar{r}^k_z).
\end{aligned}
\end{equation}
The usual approach in the interior-point literature would be to eliminate $\Delta x^k$ to obtain the so-called \textit{normal equations} since there is usually no regularization from the proximal terms, i.e. $\delta^k=0$ \cite{domahidi2012, lowenstein2024}. However, this destroys the sparsity pattern assuming our more general structure, that is, $\Phi^{-1}$ is dense due to the coupling between stages in the cost and inequalities. Instead, we eliminate $\Delta y^k$, yielding the system
\begin{equation}
\Psi \Delta x^k = \bar{r}^k
\end{equation}
where
\begin{equation}
\begin{aligned}
\Psi &= (\Phi + \delta_k^{-1} A^\top A), \\
\bar{r}^k &= \bar{r}^k_x + \delta_k^{-1} A^\top r^k_y, \\
\Delta y^k &= \delta_k^{-1}(A \Delta x_k - r^k_y). \\
\end{aligned}
\end{equation}
While this approach shares similarities with \cite{frey2020} and \cite{shin2024}, it differs fundamentally in its treatment of regularization terms. These works introduce the regularization terms heuristically, compromising either convergence guarantees or computational efficiency through required post-processing steps such as iterative refinement.

\subsection{Block-tri-diagonal-arrow Cholesky Factorization}

The symmetric matrix $\Psi$ is positive-definite and has the block-tri-diagonal-arrow structure
\begin{equation*}
\Psi \!\coloneqq\! \begin{bmatrix}
\Psi_{0,0} &\!\!\! \Psi_{1,0}^\top &\!\!\! 0 &\!\!\! \cdots &\!\!\! \Psi_{N+1,0}^\top \\[1ex]
\Psi_{1,0} &\!\!\! \Psi_{1,1} &\!\!\! \Psi_{1,2}^\top &\!\!\! \ddots &\!\!\! \Psi_{N+1,1}^\top \\
0 &\!\!\! \Psi_{1,2} &\!\!\! \ddots &\!\!\! \ddots &\!\!\! \vdots \\
\vdots &\!\!\! \ddots &\!\!\! \ddots &\!\!\! \Psi_{N,N} &\!\!\! \Psi_{N+1,N}^\top \\[1ex]
\Psi_{N+1,0} &\!\!\! \Psi_{N+1,1} &\!\!\! \cdots &\!\!\! \Psi_{N+1,N} &\!\!\! \Psi_{N+1,N+1}
\end{bmatrix}
\end{equation*}
\begin{algorithm}[t]
\caption{Construction of $\Psi$} \label{alg:psi}
\begin{algorithmic}[1]
\State $\Psi_{0,0} \leftarrow Q_0 + \rho_kI_{n_0} + \delta_k^{-1}A_{0}^\top A_{0} + C_0^\top(W^k_0+\delta_k)^{-1}C_0$
\State $\Psi_{N+1,0} \leftarrow T_0 + \delta^{-1}_k E_0^\top A_0 + F_0^\top(W^k_0+\delta_k)^{-1}C_0$
\For{$i = 0,\ldots,N-1$}
    \State $\Psi_{i+1,i+1} \leftarrow Q_{i+1} + \rho_kI_{n_{i+1}}$
    \Statex $\hspace{9.07em} + \delta_k^{-1}A_{i+1}^\top A_{i+1} + \delta_k^{-1}B_{i}^\top B_{i}$
    \Statex $\hspace{9.07em} + C_{i+1}^\top(W^k_{i+1}+\delta_k)^{-1}C_{i+1}$
    \Statex $\hspace{9.07em} + D_{i}^\top(W^k_{i}+\delta_k)^{-1}D_{i}$
    \State $\Psi_{i+1,i} \leftarrow S_i + \delta_k^{-1}B_i^\top A_i + D_i^\top (W^k_{i}+\delta_k)^{-1} C_i$
    \State $\Psi_{N+1,i+1} \leftarrow T_{i+1} + \delta_k^{-1} E_{i+1}^\top A_{i+1} + \delta_k^{-1} E_{i}^\top B_{i}$
    \Statex $\hspace{9.29em} + F_{i+1}^\top (W^k_{i+1}+\delta_k)^{-1}C_{i+1}$
    \Statex $\hspace{9.29em} + F_{i}^\top (W^k_{i}+\delta_k)^{-1}D_{i}$
\EndFor
\State $\Psi_{N+1,N+1} \leftarrow Q_g + \rho_kI_{n_g} + \sum_{j=0}^{N} \delta_k^{-1}E_j^\top E_j$
\Statex $\hspace{7.55em} + \sum_{j=0}^{N} F_j^\top(W^k_j+\delta_k)^{-1}F_j$
\end{algorithmic}
\end{algorithm}
The constructions of $\Psi$ and $\bar{r}^k$ are described in Algorithm~\ref{alg:psi} and Algorithm~\ref{alg:r_bar}, respectively.
\begin{algorithm}[t]
\caption{Construction of $\bar{r}^k$} \label{alg:r_bar}
\begin{algorithmic}[1]
\State $\bar{r}^k_0 \leftarrow r^k_{x,0} + \delta_k^{-1} A_0^\top r^k_{y,0} + C_0^\top (W_0^k+\delta_kI)^{-1} \bar{r}^k_{z,0}$
\For{$i = 1,\ldots,N$}
    \State $\bar{r}^k_{i} \leftarrow r^k_{x,i} + \delta_k^{-1} A_{i}^\top r^k_{y,i} + \delta_k^{-1} B_{i-1}^\top r^k_{y,i-1}$
    \Statex $\hspace{5.74em} + C_{i}^\top (W_{i}^k+\delta_kI)^{-1} \bar{r}^k_{z,i}$
    \Statex $\hspace{5.74em} + D_{i-1}^\top (W_{i-1}^k+\delta_kI)^{-1} \bar{r}^k_{z,i-1}$
\EndFor
\State $\bar{r}^k_{N+1} \leftarrow \sum_{i=0}^{N} \left( \delta_k^{-1} E_i^\top r^k_{y,i} + F_i^\top (W_i^k+\delta_kI)^{-1} \bar{r}^k_{z,i} \right)$
\end{algorithmic}
\end{algorithm}

The Cholesky factorization $LL^\top = \Psi$ preserves the block-tri-diagonal-arrow structure of $\Psi$, with the lower triangular factor $L$ taking the form
\begin{equation}
L \!:=\! \begin{bmatrix}
L_{0,0} &\!\!\! 0 &\!\!\! 0 &\!\!\! \cdots &\!\!\! 0 \\
L_{1,0} &\!\!\! L_{1,1} &\!\!\! 0 &\!\!\! \ddots &\!\!\! 0 \\
0 &\!\!\! L_{1,2} &\!\!\! \ddots &\!\!\! \ddots &\!\!\! \vdots \\
\vdots &\!\!\! \ddots &\!\!\! \ddots &\!\!\! L_{N,N} &\!\!\! 0 \\[1ex]
L_{N+1,0} &\!\!\! L_{N+1,1} &\!\!\! \cdots &\!\!\! L_{N+1,N} &\!\!\! L_{N+1,N+1}
\end{bmatrix}
\end{equation}
The computation of $L$ follows the procedure detailed in Algorithm~\ref{alg:factorization}.
\begin{algorithm}[t]
\caption{Factorization of $\Psi$} \label{alg:factorization}
\begin{algorithmic}[1]
\State $L_{0,0} \leftarrow \chol\left(\Psi_{0,0}\right)$
\State $L_{N+1,0} \leftarrow \Psi_{N+1,0} L_{0,0}^{-\top}$
\For{$i = 1,\ldots,N$}
    \State $L_{i,i} \leftarrow \chol\left(\Psi_{i,i} - L_{i,i-1} L_{i,i-1}^\top\right)$
    \State $L_{i,i-1} \leftarrow \Psi_{i,i-1} L_{i-1,i-1}^{-\top}$
    \State $L_{N+1,i} \leftarrow \left(\Psi_{N+1,i} - L_{N+1,i-1} L_{i,i-1}^\top\right)L_{i,i}^{-\top}$
\EndFor
\State $L_{N+1,N+1} \!\leftarrow\! \chol\left( \Psi_{N+1,N+1} \!-\! \sum_{i=0}^N L_{N+1,i} L_{N+1,i}^\top \right)$
\end{algorithmic}
\end{algorithm}
$\Delta x^k$ is computed through standard forward and backward substitution steps, as described in Algorithm~\ref{alg:solve}.

While the literature extensively documents Cholesky factorizations for both block-arrow and block-tri-diagonal matrices, the factorization of matrices with a combined block-tri-diagonal-arrow structure appears to be unexplored. This paper presents, to the best of our knowledge, the first comprehensive treatment of this specialized factorization pattern.
\begin{algorithm}[t]
\caption{Solving for $\Delta x^k$} \label{alg:solve}
\begin{algorithmic}[1]
\Statex \textbf{Forward Substitution}
\State $\Delta x^k_0 \leftarrow L_{0,0}^{-1} \bar{r}^k_0$
\For{$i = 1,\ldots,N$}
    \State $\Delta x^k_i \leftarrow L_{i,i}^{-1} \left( \bar{r}^k_i - L_{i,i-1} \Delta x^k_{i-1} \right)$
\EndFor
\State $\Delta x^k_{N+1} \leftarrow L_{N+1,N+1}^{-1}\left(\bar{r}^k_{N+1} - \sum_{j=0}^N L_{N+1,j} \Delta x^k_j\right)$
\Statex \textbf{Backward Substitution}
\State $\Delta x^k_{N+1} \leftarrow L_{N+1,N+1}^{-\top} \Delta x^k_{N+1}$
\State $\Delta x^k_{N} \leftarrow L_{N,N}^{-\top} \left( \Delta x^k_{N} - L_{N+1,N} \Delta x^k_{N+1} \right)$
\For{$i = N-1,\ldots,0$}
    \State $\Delta x^k_i \leftarrow L_{i,i}^{-1}\left(\Delta x^k_i \!-\! L_{i+1,i} \Delta x^k_{i+1} \!-\! L_{N+1,i} \Delta x^k_{N+1}\right)$
\EndFor
\end{algorithmic}
\end{algorithm}

\subsection{Flop Analysis} \label{sec:flop_analysis}

To analyze the computational complexity, we evaluate the floating point operations (flops) required by our method. To be able to compare with existing literature, we base our analysis on a variant of the extended linear quadratic control problem inspired by \cite{frison2013}
\begin{equation}
\begin{aligned}
    \min_{z, u, g} \quad & \sum_{i=0}^{N-1} \ell_i(z_i,u_i,g) + \ell_N(z_i,g) \\
    \text {s.t.}\quad & z_{i+1} = \mathcal{A}_i z_i + \mathcal{B}_i u_i + E_ig + b_i,
\end{aligned}
\end{equation}
with state $z_i \in \mathbb{R}^{n_x}$, inputs $u_i \in \mathbb{R}^{n_u}$, and system dynamics $\mathcal{A}_i \in \mathbb{R}^{n_x \times n_x}$, $\mathcal{B}_i \in \mathbb{R}^{n_x \times n_u}$. For our formulation \eqref{eq:stage_problem}, we define $x_i=(z_i,u_i)$, $A_i=\begin{bmatrix}\mathcal{A}_i & \mathcal{B}_i\end{bmatrix}$, $B_i=\begin{bmatrix}-I_{n_x} & 0_{n_x \times n_u}\end{bmatrix}$, and $S_i=0$ for $i=0,\dots,N-1$. Note, the structure of $B_i$, specifically its zero block, reduces the dimensions of the coupling matrices to $\Psi_{i+1,i}, L_{i+1,i} \in \mathbb{R}^{n_x \times (n_x+n_u)}$ rather than $\mathbb{R}^{(n_x+n_u) \times (n_x+n_u)}$.

Using the elementary matrix operation costs from Table~\ref{table:elem_flops}, we derive the computational complexity of both factorization and solution steps, presented in Table~\ref{table:fact_solv_flops}. For $n_g=0$, our results align with previous findings \cite{malyshev2018}, and notably, the factorization exhibits identical complexity as the Riccati recursion implemented in HPIPM \cite{frison2020}. The computationally intensive construction of $\Psi$ can be optimized by caching invariant terms during the initial setup phase.
\begin{table}
\centering
\caption{Cost of elementary matrix operations with matrices\\
$A : m \times n, B : n \times p, D, L : n \times n, L$ lower triangular.}
\vspace{-1.3em}
\begin{tabular}{lll}
\toprule
Operation & Kernel & Cost (flops) \\
\midrule
Matrix matrix multiplication $A \cdot B$ & \texttt{gemm} & $2mnp$ \\
Symmetric matrix multiplication $A \cdot A^T$ & \texttt{syrk} & $m^2n$ \\
Cholesky decomposition s.t. $LL^T = D$ & \texttt{potrf} & $1/3n^3$ \\
Solving triangular matrix $A \cdot L^{-T}$ & \texttt{trsm} & $mn^2$ \\
\bottomrule
\end{tabular}
\label{table:elem_flops}
\end{table}

\begin{table}
\centering
\caption{Cost of factorization and solve for MPC problems}
\vspace{-1em}
\begin{tabular}{lll}
\toprule
Operation & Cost (flops) \\
\midrule
Construction of $\Psi$ & $N(4n_x^3+4n_x^2n_u+n_xn_u^2)$ \\
& $+Nn_g(n_gn_x+n_x^2+n_xn_u)$ \\
$^\dagger$Factorization of $\Psi$ & $N(\frac{7}{3}n_x^3+4n_x^2n_u+2n_xn_u^2+\frac{1}{3}n_u^3)$ \\
& $+Nn_g(3n_x^2+4n_xn_u+n_u^2) + \frac{1}{3}n_g^3$ \\
\midrule
Construction of $\bar{r}^k$ & $N(2n_x^2+n_xn_u+n_xn_g)$ \\
Solving for $\Delta x^k$ & $N(4n_x^2+6n_xn_u+2n_u^2)$ \\
& $+2Nn_g(n_x+n_u)+2n_g^2$ \\
Reconstruction of $\Delta y^k$ & $N(2n_x^2+n_xn_u+n_xn_g)$ \\
\bottomrule
\end{tabular}
\label{table:fact_solv_flops}
\vspace{-1em}
\end{table}

Existing OCP structure-exploiting algorithms, such as HPIPM \cite{frison2020}, lack direct support for the global decision variable $g$ in \eqref{eq:stage_problem}. To accommodate this limitation, the conventional approach embeds $g$ into an augmented stage state $x_i=(z_i,u_i,g)$, reformulating the dynamics matrices as
\begin{equation*}
\begin{aligned}
    A_i&=\begin{bmatrix}
    \mathcal{A}_i & \mathcal{B}_i & E_i \\
    0_{n_g \times n_x} & 0_{n_g \times n_u} & I_{n_g}
    \end{bmatrix}, \\
    \quad B_i&=\begin{bmatrix}
    -I_{n_x} & 0_{n_x \times n_u} & 0_{n_x \times n_g} \\
    0_{n_g \times n_x} & 0_{n_g \times n_u} & -I_{n_g}
    \end{bmatrix}.
\end{aligned}
\end{equation*}
Thus, subsituting $n_x \leftarrow n_x+n_g$ and $n_g \leftarrow 0$ in Table~\ref{table:fact_solv_flops}~($\dagger$) results in an additional cost of
\begin{equation*}
    Nn_g(\frac{7}{3}n_g^2 + 7n_gn_x + 4n_gn_u + 4n_x^2 + 4n_xn_u + n_u^2) - n_g^3
\end{equation*}
flops compared to the factorization of $\Psi$ before augmentation. This not only incurs a strictly higher factorization cost, but also increases the $\mathcal{O}(n_g^3)$ term in Table~\ref{table:fact_solv_flops} to $\mathcal{O}(Nn_g^3)$.

\section{Numerical Implementation} \label{sec:numerical_implementation}

We have implemented the factorization routines described in Section~\ref{sec:solving_kkt} as a new backend within the solver PIQP \cite{schwan2023}. The implementation utilizes C/C++ with the Eigen3 library for the default sparse backend, while leveraging BLASFEO \cite{frison2018} for optimized linear algebra operations in the new backend. Consistent with the original solver, we maintain Python, Matlab, and R interfaces.

\subsection{Automatic Structure Detection}

\begin{revenv}    
To enhance usability and facilitate adoption, our implementation preserves compatibility with the generic sparse interface of problem \eqref{eq:general_qp} through an automatic structure detection algorithm. The detection process begins by constructing the KKT matrix sparsity pattern $\Psi$ from the problem matrices as
\begin{equation}
\Psi = P + I + A^\top A + G^\top G,
\end{equation}
where we only retain the lower triangular structure to identify the underlying block arrow pattern.

The algorithm employs a greedy approach that iterates through each row of $\Psi$, making local decisions about block partitioning based on estimated factorization costs. For each potential block configuration, we estimate the computational cost of tridiagonal factorization (diagonal blocks) versus arrow factorization (corner block), including operations such as Cholesky decomposition (\texttt{potrf}), matrix multiplication (\texttt{gemm}), triangular solve (\texttt{trsm}), and symmetric rank-$k$ update (\texttt{syrk}). The algorithm maintains a running estimate of the current block structure and greedily assigns matrix elements to either extend the current diagonal block or increase the arrow width based on which choice minimizes the total FLOP count. The complete procedure is outlined in Algorithm~\ref{alg:structure_detection}.

\begin{algorithm}[t]
\begin{revenv}
\caption{Automatic Arrow Structure Detection} \label{alg:structure_detection}
\begin{algorithmic}[1]
\State Construct sparsity pattern $\Psi \leftarrow P + I + A^\top A + G^\top G$
\State Initialize block info: $\text{start} \leftarrow 0$, $\text{diag\_size} \leftarrow 0$, $\text{off\_diag\_size} \leftarrow 0$, $\text{arrow\_width} \leftarrow 0$
\For{$i = 0,\ldots,n-1$}
    \For{each nonzero $(i,j)$ in row $i$ of $\Psi$}
        \State Compute FLOP cost for extending diagonal block
        \State Compute FLOP cost for increasing arrow width  
        \If{diagonal cost $\leq$ arrow width cost}
            \State Update diagonal and off-diagonal block sizes
        \Else
            \State Update arrow width
        \EndIf
    \EndFor
    \If{block completion criteria met}
        \State Finalize current block and add to block structure
        \State Initialize next block with updated parameters
    \EndIf
\EndFor
\State Add final arrow block to structure
\end{algorithmic}
\end{revenv}
\end{algorithm}

Block completion occurs when the algorithm reaches the arrow region or the next iteration would significantly increase the block size. After initial block detection, the algorithm performs a merging step to consolidate adjacent blocks that can be efficiently combined without increasing the total FLOP count.

\begin{figure}[tbp!]
\centering
\includegraphics[width=0.3\textwidth]{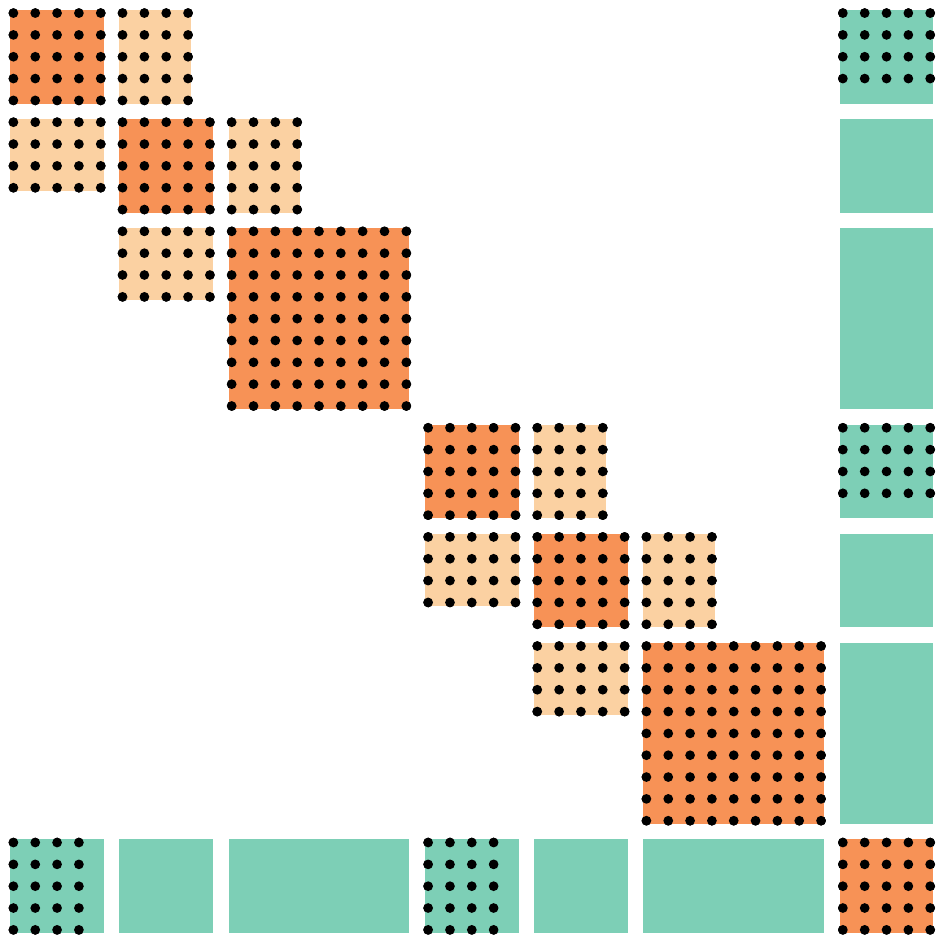}
\caption{Sparsity structure and detected blocks of $\Psi$ for scenario problem \eqref{eq:scenario_problem} with $M=2$, $N=4$, and $N_s=2$, \rev{where dots are non-zero entires.}}
\label{fig:arrow_sparsity}
\end{figure}

Figure~\ref{fig:arrow_sparsity} illustrates the resulting block sparsity detection for the scenario problem \eqref{eq:scenario_problem}. Importantly, while the sparsity pattern of $\Psi$ depends on the ordering of decision variables, it remains invariant to constraint ordering due to the quadratic terms $A^\top A$ and $G^\top G$ that are invariant to row permutations. This characteristic significantly simplifies the user experience, requiring only correct decision variable ordering.

This greedy approach always finds a feasible block structure but may not be globally optimal, as the underlying variable ordering problem is known to be NP-complete \cite{Yannakakis1981CtM}. However, the local cost-based decisions provide effective heuristics for practical problems with natural block structure.
\end{revenv}

\section{Numerical Examples} \label{sec:numerical_examples}

We benchmark our algorithm against the original PIQP (sparse) \cite{schwan2023} solver, the generic QP solvers OSQP \cite{stellato2020} and QPALM \cite{hermans2019}, and the tailored OCP solver HPIPM \cite{frison2020}. To assess the impact of hardware optimization, we evaluate PIQP using both SSE and AVX2 instruction sets. HPIPM is compiled targeting the AVX2 instruction set \rev{on an x86 architecture}. All benchmarks are executed on an Intel Core i9 2.4 GHz CPU with Turbo Boost disabled to ensure consistent thermal performance.

Following established literature, we base our numerical experiments on the classical spring-mass system \cite{wang2010}, extended with an additional cost term penalizing input differences between consecutive stages. The resulting optimization problem is formulated as
\begin{equation} \label{eq:chain_mass_problem}
\begin{aligned}
    \min_{z, u} \quad & \sum_{i=0}^{N-1} \ell_i(z_i,u_i,u_{i+1}) + \ell_N(z_N) \\
    \text {s.t.}\quad & z_0 = \bar{x}_0, \\
    & z_{i+1} = A z_i + B u_i, \\
    & -4 \cdot \mathbf{1}_{n_x} \leq z_i \leq 4 \cdot \mathbf{1}_{n_x}, \\
    & -0.5 \cdot \mathbf{1}_{n_u} \leq u_i \leq 0.5 \cdot \mathbf{1}_{n_u},
\end{aligned}
\end{equation}
with stage cost
\begin{equation*}
\begin{aligned}
    \ell_i(z_i, u_i, u_{i+1}) &= z_i^\top Q z_i + u_i^\top R u_i \\
    &+ \left( u_i - u_{i+1} \right)^\top R_d \left( u_i - u_{i+1} \right),
\end{aligned}
\end{equation*}
and terminal cost $\ell_N(z_N) = z_N^\top Q_N z_N$, where $z_i \in \mathbb{R}^{n_x}$ and $u_i \in \mathbb{R}^{n_u}$. For a system with $M$ masses, the state dimension is $n_x=2M$, with $M-1$ actuators yielding an input dimension of $n_u=M-1$. We select cost matrices $Q=10^3 \cdot I_{n_x}$, $R = R_d=10^{-1} \cdot I_{n_u}$, and derive $Q_N$ as the solution from the discrete-time Riccati equation. To exploit the block-tri-diagonal structure in $\Psi$, we organize the decision variables as $x = (z_0,u_0,\dots,z_{N-1},u_{N-1},z_N)$ with no global variables, i.e. $g$ empty.

Following \cite{lowenstein2024}, we initialize the system state by sampling $\gamma \sim U[0.5,1.5]$ and subsequently drawing $\bar{x}_0 \sim U[-\gamma, \gamma]$, ensuring a balanced distribution of active and inactive constraints. All solvers are configured with identical absolute and relative tolerances of $10^{-6}$ and warm-started using zero-input trajectories simulated from $\bar{x}_0$.

\begin{figure}[tbp!]
\centering
\includegraphics[width=0.5\textwidth]{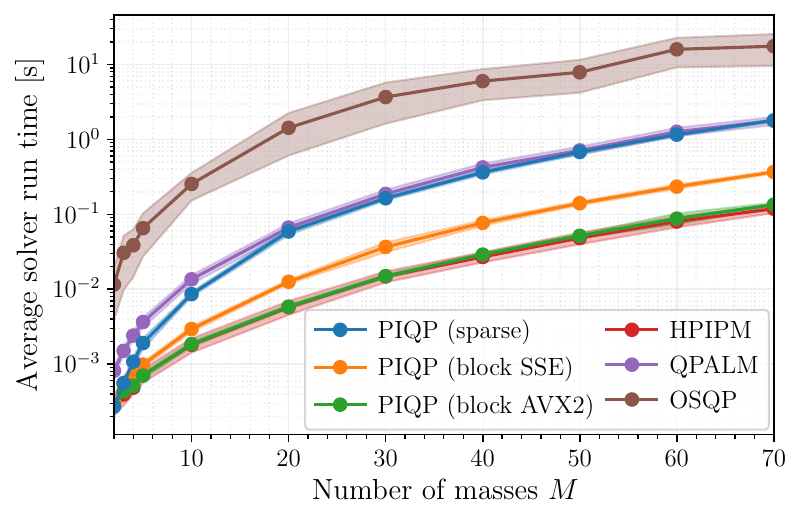}
\vspace{-2em}
\caption{Average solve time and standard dedication over 30 runs for a varying number of masses $M$ for $R_d=0$ and $N=15$.}
\label{fig:benchmark_M2-70_N15_default}
\end{figure}

Our initial benchmark considers the case where $R_d=0$, eliminating input coupling between stages. For each configuration, we perform 30 runs with the number of masses $M$ ranging from 2 to 70. The results in Figure~\ref{fig:benchmark_M2-70_N15_default} demonstrate that our new PIQP backend achieves substantial performance improvements: a 13x speed-up compared to the generic sparse backend and up to 245x acceleration versus OSQP. The transition from SSE to AVX2 instructions yields the theoretical 2x performance gain, highlighting the effectiveness of dense kernels and vectorization. For small problem sizes, the vectorization and cache locality of the dense kernels are less advantageous, i.e., the generic sparse interface might even be slightly faster. \rev{We can also see on-par performance with HPIPM, which can be expected due to HPIPM also using an interior-point method, BLASFEO as linear algebra backend, and has the same expected flop count as discussed in Section~\ref{sec:flop_analysis}.}

\begin{figure}[tbp!]
\centering
\includegraphics[width=0.5\textwidth]{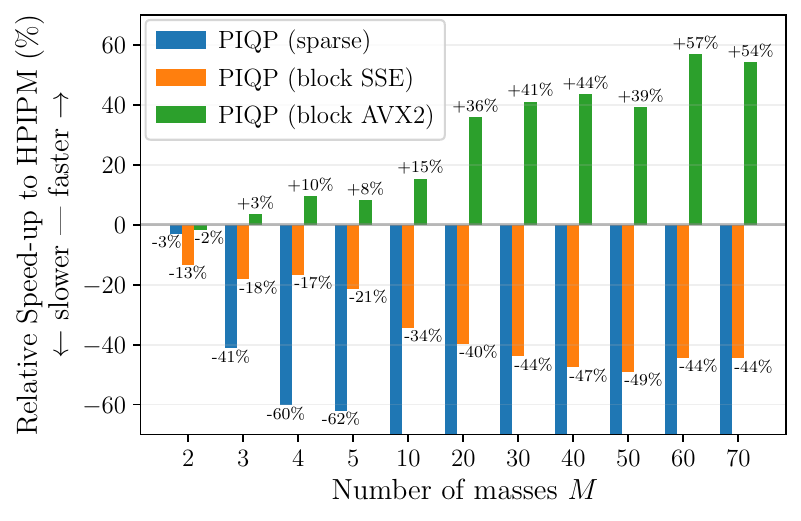}
\vspace{-2em}
\caption{Relative speed-up compared to HPIPM for a varying number of masses $M$ for $R_d=10^{-1} \cdot I_{n_u}$ and $N=15$.}
\label{fig:benchmark_M2-70_N15_cost_diff}
\end{figure}

Introducing input coupling with $R_d=10^{-1}\cdot I_{n_u}$ challenges the standard OCP structure expected by HPIPM. Similar as in Section~\ref{sec:flop_analysis}, we address this by augmenting the state to $x_i=(z_i,u_i,\bar{u}_{i+1})$ with the constraint $\bar{u}_i = u_i$, effectively decoupling the stages in the cost. Under these conditions, as shown in Figure~\ref{fig:benchmark_M2-70_N15_cost_diff}, PIQP outperforms HPIPM by up to 57\%. \rev{The slightly worse performance for smaller problem instances is mainly due to overhead in the implementation of PIQP, introduced by the more general interface and support for general QP problems.}

To evaluate performance with global decision variables $g$, we extend problem \eqref{eq:chain_mass_problem} to a robust scenario formulation \cite{lucia2013}
\begin{equation} \label{eq:scenario_problem}
\begin{aligned}
    \min_{z, u} \quad & \frac{1}{N_s} \sum_{j=1}^{N_s} \left( \sum_{i=0}^{N-1} \ell_i(z_i^j,u_i^j) + \ell_N(z_N^j) \right) \\
    \text {s.t.}\quad & z_0^j = \bar{x}_0, u_0^1 = u_0^j \\
    & z^j_{i+1} = A^j z_i^j + B^j u_i^j, \\
    & -4 \cdot \mathbf{1}_{n_x} \leq z^j_i \leq 4 \cdot \mathbf{1}_{n_x}, \\
    & -0.5 \cdot \mathbf{1}_{n_u} \leq u^j_i \leq 0.5 \cdot \mathbf{1}_{n_u},
\end{aligned}
\end{equation}
with $N_s$ scenarios. The system matrices $A^j$ and $B^j$ describe $N_s$ distinct spring-mass systems, each with spring constants sampled uniformly from $U[1,2]$. Organizing the variables as $x=(z_1^1,u_1^1,\dots,z_{N}^1,z_1^2,u_1^2,\dots,z_{N}^{N_s})$ and $g=(x_0,u_0)$ yields the block-tri-diagonal-arrow structure illustrated in Figure~\ref{fig:arrow_sparsity}. The performance comparison between the block-sparse backend using AVX2 instructions and the sparse backend reveals increasing efficiency with problem size, reaching \rev{up to 7.9x overall} speed-up through improved vectorization, cache utilization, and linear memory access patterns that better utilize available memory bandwidth on modern CPUs.

\begin{figure*}[tbp!]
\centering
\includegraphics[width=1\textwidth]{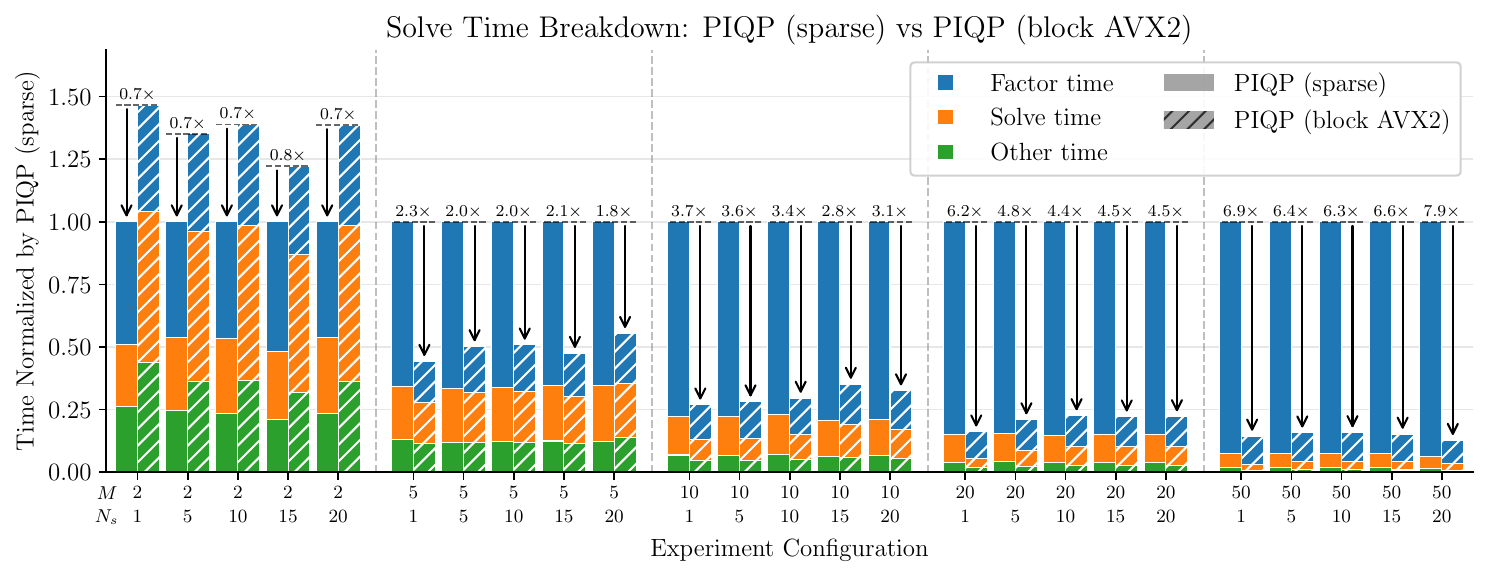}
\vspace{-2em}
\caption{Speed-up of the block-sparse backend compared to the sparse backend in PIQP for different number of masses $M$, scenarios $N_s$ and horizon $N=15$. \rev{The stacked bars show the time breakdown for each solver, normalized with respect to the sparse backend. Blue represents KKT factorization time, orange the KKT solve time, and green all remaining operations. The waterfall arrows indicate the overall speed-up factor.}}
\label{fig:scenario_speedup_heatmap}
\end{figure*}

\begin{revenv}
Figure~\ref{fig:scenario_speedup_heatmap} presents a comprehensive breakdown of the performance gains achieved by the block-sparse backend. Beyond the overall end-to-end speedup, we provide detailed performance metrics for the factorization time and solve routines of the KKT system, as well as the remaining computations (e.g., line search, residual calculation, etc.). The results reveal a clear performance scaling pattern across problem sizes, with the block-sparse implementation achieving increasingly substantial speedups as problem dimensions grow. For smaller problems, the sparse backend maintains a slight advantage, mainly due to less overhead and lower solve times. However, as the problem size increases, the block-sparse backend outperforms the sparse backend with overall speedups of up to 7.9x.

The breakdown shows that while factorization time reductions contribute significantly to overall performance gains, improvements in solve time and other computational components provide less and might actually be slower for small problem instances due to overhead.
\end{revenv}

\bibliographystyle{IEEEtran}
\bibliography{refs}

\end{document}